\documentclass[11pt,3p, authoryear]{elsarticle}
\usepackage{graphicx}

\usepackage{algorithm}  
\usepackage{algorithmic}
\usepackage{amssymb,url,amsmath}
\usepackage{bm} 

\usepackage[dvipsnames]{xcolor}
\usepackage{booktabs}	
\usepackage{epstopdf}
\usepackage{color}
\usepackage{pifont}
\usepackage{tikz}
\usepackage{hyperref}
\usepackage{eurosym}

\usepackage{pgfplots}
\usepackage{setspace}
\usepackage{rotating}
\usepackage{circuitikz}
\usepackage{bbm}

\newtheorem{theorem}{Theorem}
\newtheorem{inequality}{Inequality}

\newproof{pf}{Proof}
\tikzset{>=latex}

\usepackage{bm}
\newcommand*{\QEDA}{\hfill\ensuremath{\blacksquare}}

\newcommand{\CC}{\mathcal{C}}

\renewcommand{\CC}{\mathcal{C}}
\newcommand{\D}{\mathcal{D}}
\newcommand{\HC}{\mathcal{H}}

\allowdisplaybreaks
\makeatother
\makeatletter
\def\ps@pprintTitle{%
 \let\@oddhead\@empty
 \let\@evenhead\@empty
 \def\@oddfoot{}%
 \let\@evenfoot\@oddfoot}
\makeatother

\begin{document}
\begin{frontmatter}

\title{A comment on ``Performance guarantees of a greedy algorithm for minimizing a supermodular set function on comatroid"}

\author[a1]{Orcun Karaca}
\ead{okaraca@ethz.ch}
\author[a2]{Baiwei Guo}
\ead{baiwei.guo@epfl.ch}
\author[a1]{Maryam Kamgarpour}
\ead{maryamk@ethz.ch}
\address[a1]{Automatic Control Laboratory, D-ITET, ETH Z{\"u}rich, Switzerland}
\address[a2]{Automatic Control Laboratory, EPFL, Switzerland.}
\vspace{-.1cm}
\begin{abstract}
We provide a counterexample to the performance guarantee obtained in the paper ``Performance guarantees of a greedy algorithm for minimizing a supermodular set function on comatroid", which was published in Volume 171 of the European Journal of Operational Research. We comment on where this error originates from in the proof of the main theorem.
\end{abstract}
\begin{keyword}
Combinatorial optimization \sep Greedy algorithm \sep Matroid theory
\end{keyword}
\end{frontmatter}
\section{Problem formulation}
Let $U$ be a finite set with $|U|=n$. Let hereditary system $\HC =(U,\D)$ denote a comatroid, where the family $\D$ denotes the dependence system, and let $f:2^U\rightarrow \mathbb{R}_+$ denote a set function.
\citet{il2006performance} consider solving 
\begin{equation}\label{eq:main}\min\{f(X):X\in\CC\},\end{equation}
where $\CC$ is the family given by all circuits of a comatroid $\HC=(U,\D)$ with girth $p$, and $f$ is supermodular, nonincreasing, and $f(U)=0.$ We refer to \citep{il2006performance} and \citep{il2003hereditary} for the definitions relating to comatroids and hereditary systems. 

This problem is known to be NP-hard since the well-known $p$-median problem can be captured as a special case.
As a heuristic, \citet{il2006performance} propose the greedy descent algorithm (also known as reverse greedy or stingy algorithm), which proceeds as follows:\\

\noindent\textit{Greedy descent algorithm:}\\
\textit{Step $0$:} Set $X_0 = U.$ Go to Step 1.\\
\textit{Step $i$:} ($i\geq1$): Select $x_i\in X_{i-1}$ such that
    \begin{equation}\label{eq:exploit}d_{x_i}(X_{i-1})=\min_{\substack{x\in X_{i-1}\\ X_{i-1}\setminus \{x\} \in \D}} d_x(X_{i-1}),\end{equation}
    where $d_x(X)=f(X\setminus\{x\}) - f(X)$ is the marginal increment in $f$ when removing $\{x\}$ from the set~$X$. Set $X_i = X_{i-1}\setminus \{x_i\}$. If $i=n-p$, then stop. Otherwise go to Step $i+1.$\\
    \textit{End.}\\

The paper contained the following theorem regarding the suboptimality bound of the greedy heuristic when applied to \eqref{eq:main}.
\begin{theorem}\label{thm:main}\citep[Theorem~1]{il2006performance}
Let OPT be an optimal solution to~\eqref{eq:main} on an arbitrary comatroid and GR be the solution returned by the greedy descent algorithm. Then,
$$\dfrac{f(\text{GR})}{f(\text{OPT})}\leq \dfrac{1}{t} \left(\left(1+\dfrac{t}{q}\right)^q-1\right),$$
where $q=n-p$, $t=s/1-s$, and $s$ is the solution to the following problem
\begin{equation*}
    s = \max_{\substack{x\in U,\\ f(\{x\})< f(\emptyset)}} \dfrac{(f(\emptyset)-f(\{x\}))-(f(U\setminus \{x\})-f(U))}{(f(\emptyset)-f(\{x\}))}.
\end{equation*}
\end{theorem}

In the following, we provide a counterexample showing that the guarantee in Theorem~\ref{thm:main} does not necessarily hold. Then, we comment on the mistake found in \citep[Lemma~1]{il2006performance}, which is then utilized in constructing the linear program \citep[equation~(8)]{il2006performance} in the proof of Theorem~\ref{thm:main}.
\section{Counterexample}
Set $n=4$, $U=\{1,2,3,4\}$. Consider the following nonincreasing supermodular function:
\begin{equation*}
    \begin{split}
        &f(\emptyset)= 6,\ f(\{1\})= 4,\ f(\{2\})= 5,\ f(\{3\})= 4,\ f(\{4\})= 4,\\
        &f(\{3,4\})= 2,\ f(\{1,2\})= 3,\ f(\{2,4\})= 3,\ f(\{2,3\})= 3,\ f(\{1,4\})= 2,\ f(\{1,3\})= 2,\\
        &f(\{1,2,3\})= 1,\ f(\{2,3,4\})= 1,\ f(\{1,3,4\})= 1,\ f(\{1,2,4\})= 1, \\
        & f(\{1,2,3,4\}) = f(U) = 0.
    \end{split}
\end{equation*}
Compute the steepness $s$ of function~$f$:
$$s=\dfrac{(6-4)-(1-0)}{(6-4)}=0.5.$$
Hence, $t=1.$
Define the comatroid $\HC=(U,\D)$ as follows:
\begin{equation*}
    \begin{split}
       \D =&\{ U=\{1,2,3,4\},\\
       &\{2,3,4\},\{1,3,4\},\{1,2,4\},\{1,2,3\},\\
       &\{1,2\},\{1,4\},\{2,3\},\{3,4\}\}.
    \end{split}
\end{equation*}
Note that girth is given by $p=2$, thus $q=2$. This comatroid was previously studied in \citep[Remark~3]{il2006performance}. Clearly, the family given by all circuits of this comatroid is $\CC=\{\{1,2\},\{1,4\},\{2,3\},\{3,4\}\}.$

Consider \eqref{eq:main} with the comatroid $\HC=(U,\D)$ and the objective $f$. The greedy descent algorithm can find the solution $\text{GR}=\{1,2\}$\footnote{In this example, the selection done in Step 1 by \eqref{eq:exploit} is not unique. We can obtain $\text{GR}=\{1,2\}$ (or $\text{GR}=\{2,3\}$) if $X_{1}=\{1,2,3\}$ is chosen in Step 1. Notice that if $X_{1}=\{1,3,4\}$, greedy descent can find an optimal solution. This statement also holds for $X_{1}=\{2,3,4\}$ and $X_{1}=\{1,2,4\}$.}, whereas an optimal solution is given by $\text{OPT}=\{3,4\}$. Theorem~\ref{thm:main} claims
\begin{equation*}
\dfrac{3}{2}=\dfrac{f(\text{GR})}{f(\text{OPT})}  \leq  \dfrac{1}{t} \left(\left(1+\dfrac{t}{q}\right)^q-1\right) = \dfrac{1}{1} \left(\left(1+\dfrac{1}{2}\right)^2-1\right)= 1.25,
\end{equation*}
which is not correct.
\section{The error in the proof of the main theorem}
Denote the complements by $\overline{X}=U\setminus X$. The proof of \citep[Theorem~1]{il2006performance} relies on 
\citep[Lemma~1]{il2006performance}. This lemma exploits the following inequality in its proof.

\begin{inequality}\label{fact:d}
\begin{equation*}
\sum_{b\in \overline{\text{OPT}}\setminus\overline{X}_{i-1}}d_{b}(X_{i-1})\geq |\overline{\text{OPT}}\setminus\overline{X}_{i-1}| d_{x_i}(X_{i-1}).\end{equation*}
\end{inequality}

The above inequality and the resulting \citep[Lemma~1]{il2006performance} is then utilized in constructing the linear program in \citep[equation~(8)]{il2006performance}. However, Inequality~\ref{fact:d} is not necessarily true. For instance, in the above counterexample, $x_2=3$, $X_{1}=\{1,2,3\}$, $\text{OPT}=\{3,4\}$, and $\overline{\text{OPT}}\setminus\overline{X}_{1}=\{1,2\}.$ For $i=2$, inserting  $d_{1}(X_{1})=2$, $d_{2}(X_{1})=1$, and $d_{3}(X_{1})=2$ into Inequality~\ref{fact:d}, we obtain
\begin{equation*}
d_{1}(X_{1})+d_{2}(X_{1})=2+1\geq 2\times 2=2\times d_{3}(X_{1}),\end{equation*}
which is not correct.

As an insight, the authors conclude Inequality~\ref{fact:d} using the following statement:

``By~\eqref{eq:exploit},  for every $b\in \overline{\text{OPT}}\setminus\overline{X}_{i-1}$, it holds that $d_b(X_{i-1})\geq d_{x_i}(X_{i-1})."$

This is equivalent to
$$d_b(X_{i-1})\geq d_{x_i}(X_{i-1}),\ \forall b\in {X}_{i-1}\setminus \text{OPT}.$$
For the above inequality to hold by~\eqref{eq:exploit}, for any $b\in {X}_{i-1}\setminus \text{OPT}$, we should have ${X}_{i-1}\setminus b\in \D.$ This is not necessarily true. For instance, in the above example, $X_{1}=\{1,2,3\}$ and $\text{OPT}=\{3,4\}$, but $2\in {X}_{i-1}\setminus \text{OPT}$ is not considered by the step found in~\eqref{eq:exploit}, since $X_{1}\setminus\{2\}\notin \D.$

\section{A correction to the error in Inequality~\ref{fact:d}}

In our work~\citep{guo2019actuator}, we studied a greedy heuristic for a problem similar to~\eqref{eq:main} where the constraint set is instead the base of a matroid, and the objective is neither supermodular nor submodular but characterized by the notions of curvature and submodularity ratio. Invoking ideas from \citep[Lemma~4]{guo2019actuator}, it is possible to revise and correct Inequality~\ref{fact:d}.
\begin{inequality}\label{fact:new}
$$\sum_{b\in \overline{\text{OPT}}\setminus\overline{X}_{i-1}}d_{b}(X_{i-1})\geq (q-(i-1)) d_{x_i}(X_{i-1}).$$
\end{inequality}
\begin{pf}
\setstretch{1.25}
From the properties of comatroids derived originally in \citep[Theorem~2: Statement~(D2)]{il2003hereditary}, it can be verified that there exist $|X_{i-1}|-|\text{OPT}|$ distinct elements from ${X}_{i-1}\setminus \text{OPT}$ such that after the exclusion of these elements from ${X}_{i-1}$, we still obtain a set that lies in the comatroid. Let $R\subseteq{X}_{i-1}\setminus \text{OPT}$ denote one such subset with exactly $|X_{i-1}|-|\text{OPT}|=n-(i-1) - p=q-(i-1)$ elements. We then obtain the following,
$$\sum_{b\in \overline{\text{OPT}}\setminus\overline{X}_{i-1}}d_{b}(X_{i-1})=\sum_{b\in{X}_{i-1}\setminus \text{OPT}}d_{b}(X_{i-1})\,\geq\, \sum_{b\in R}d_{b}(X_{i-1})\,\geq\, (q-(i-1)) d_{x_i}(X_{i-1}).$$
The first equality comes from $\overline{\text{OPT}}\setminus\overline{X}_{i-1}={X}_{i-1}\setminus \text{OPT}$. The first inequality follows since function~$d$ maps to nonnegative real numbers. The fact that $X_{i-1}\setminus \{x\} \in \D$ for all $x\in R$ and the definition in \eqref{eq:exploit} conclude the last inequality.\hfill\QEDA
\end{pf}

In contrast to Inequality~\ref{fact:d}, Inequality~\ref{fact:new} involves a factor that is independent of the greedy and the optimal solutions. However, this factor is smaller since $(q-(i-1))\leq |\overline{\text{OPT}}\setminus\overline{X}_{i-1}|.$ Revisiting our counterexample, for $i=2$ inserting  $d_{1}(X_{1})=2$, $d_{2}(X_{1})=1$, and $d_{3}(X_{1})=2$ into Inequality~\ref{fact:new}, we obtain
\begin{equation*}
d_{1}(X_{1})+d_{2}(X_{1})=2+1\geq 1\times 2=1\times d_{3}(X_{1}),\end{equation*}
which is correct.

As a future work, it would be interesting to refine \citep[Theorem~1]{il2006performance} according to Inequality~\ref{fact:new}.
Finally, note that \citep[Theorems 2 and 3]{il2006performance} are also wrong since they are corollaries of the results of \citep[Lemma~1 and Theorem~1]{il2006performance}.
\bibliography{library}
\end{document}